\documentclass[12pt]{article}
\usepackage{amsmath, amssymb, amsthm}

\theoremstyle{definition}
\newcommand{\fone}{{\mathbb{F}_q}}

\newcommand{\pone}{{\mathcal{P}^2(\fone)}}

\numberwithin{equation}{section}

\newtheorem{thm}{THEOREM}[section]

\newtheorem{lem}[thm]{Lemma}
\newtheorem{cor}[thm]{Corollary}
\newtheorem{prop}[thm]{PROPOSITION}

\theoremstyle{definition}
\newtheorem{defn}{Definition}[section]
\theoremstyle{remark}
\newtheorem{rem}{Remark}[section]

\newcommand{\tref}[1]{Theorem~\ref{#1}}
\newcommand{\cref}[1]{Korollar~\ref{#1}}

\newcommand{\lref}[1]{Lemma~\ref{#1}}

\title{Polyhedra with specified links.}
\author{\sc Alina Vdovina}
\date{
\small  {Mathematisches Institut}\\
\small  {Beringstrasse 1, 53115 Bonn}\\
\small  { e-mail: \ \ alina@math.uni-bonn.de}
}

\begin{document}

\maketitle

\abstract{
We construct compact polyhedra with $m$-gonal faces whose links are
generalized $3$-gons.  It gives examples of cocompact hyperbolic bildings
of type $P(m,3)$.  For $m=3$ we get compact spaces covered by Euclidean
buildings of type $A_2$.
}

\section{Introduction}

\subsection{Preliminaries} 
Given a graph $G$ we assign to each edge the length $1$. The
diameter of the graph is its diameter as a length metric space,
its  injectivity radius is half of the length of the smallest circuit.

Due to \cite{BB}, \cite{CL} or \cite{Ronan} the following 
definition is equivalent to the usual one

\begin{defn}
For a natural number $m$ we call a connected graph $G$ a generalized
$m$-gon, if its diameter and injectivity radius are both equal to $m$.
\end{defn}

A graph is {\em bipartite} if its set of vertices
can be partitioned into two disjoint subsets $P$ and $L$ such that no
two vertices in the same subset lie on a common edge.
Such a graph can be interpreted as a planar geometry, i.e.
a set of points $P$ and a set of lines $L$ and an incidence relation
$R\subset P\times L$. On the other hand each planar geometry
can be considered as a bipartite graph.

Under this correspondence projective planes are the same as
generalized $3$-gones (\cite{Ronan}).

Let $G$ be a planar geometry.
For a line $y\in L$ we denote by $I(y)$ the set of all
points $x\in P$ incident to $y$. If no confusion can arise
we shall write $x\in y$ instead of $x\in I(y)$ and $y_1\cap
y_2$ instead of $I(y_1)\cap I(y_2)$.
A subset $S$ of $P$ is called collinear if it is contained
in some set $I(y)$, i.e. if all points of $S$ are incident to a line.

Given a planar geometry $G$ we shall denote by $G^{\prime}$ its
dual geometry arising by calling lines resp. points of $G$ points
resp. lines of $G^{\prime}$. The  graphs correspomding to
$G$ and  $G^{\prime}$ are isomorphic.

We will call a {\em polyhedron} a two-dimensional
complex which is obtained from several oriented $p$-gons
by identification of corresponding sides.
Consider a point of the polyhedron and 
take a sphere of a small radius at this point.
The intersection of the sphere with the polyhedron is
a graph, which is called the {\em link} at this point.

\begin{defn}  Let $\mathcal{P}(p,m)$ be a tessellation of the
hyperbolic plane by  regular polygons with $p$ sides,
with angles $\pi/m$ in each vertex where $m$ is an integer.
A {\em hyperbolic building} of type  $\mathcal{P}(p,m)$
is a polygonal complex $X$,
which can be expressed as the union of subcomplexes called 
apartments such that:

1. Every apartment is isomorphic to $\mathcal{P}(p,m)$.

2. For any two polygons of $X$, there is an apartment
containing both of them.

3. For any two apartments $A_1, A_2 \in X$ containing
the same polygon, there exists an isomorphism $ A_1 \to A_2$
fixing $A_1 \cap A_2$.
\end{defn}

If we replace in the above definition the tessalition $\mathcal{P}(p,m)$ of
the hyperbolic plane by the tessalation $\tilde A_2$ of the Euclidean plane
by regular triangles we get the definition of the Euclidean building of type 
$A_2$.

Let $C_p$ be a polyhedron whose faces are $p$-gons
and whose links are generalized $m$-gons with $mp>2m+p$.
We equip every face of $C_p$ with the hyperbolic
metric such that all sides of the polygons are geodesics
and all angles are $\pi/m$.
Then the universal covering of such a
polyhedron is a hyperbolic building, see \cite{Paulin}.

In the case $p=3$, $m=3$, i.e. $C_p$ is a simplicial polyhedron, 
we can give a Euclidean metric to every face. In this metric
all sides of the triangles
are geodesics of the same length. The universal coverings of these polyhedra
are Euclidean buildings, see \cite{BB}, \cite{Ba}, \cite{CL}.

So, to construct hyperbolic and Euclidean buildings
with compact quotients, it is sufficient to construct
finite polyhedra with appropriate links.

The main result of the  paper is a construction of a family of 
compact polyhedra
with $m$-gonal faces (for any $m \geq 3$) whose links
are generalized $3$-gons.

One of the main tools is a bijection $T$ of a 
special type between points and lines
of a finite projective plane $G$.
If such a bijection exists, we can construct a family of compact polyhedra
with $m$-gonal faces, with any $m \geq 3$ whose links
are generalized $3$-gons. 
The existence of $T$ in known for the projective planes
over finite fields of characteristique $\neq 3$ (chapter 3).
But the projective plane of order $3$ such a bijection exists
as well.

So, if one can prove the existence of $T$ for
a finite projective plane $G$ (even nondesarguesian),
then chapters 2.2 and 2.3 immediately give the existence
of buildings with $G$ as the link.

\subsection{Polygonal presentation.}

We recall the definition of the polygonal presentation, given in \cite{V}.

\noindent
{\bf Definition.} Suppose we have $n$ disjoint connected bipartite graphs
 $G_1, G_2, \ldots G_n$.
Let $P_i$ and $L_i$ be the sets of black and white vertices respectively in
$G_i$, $i=1,...,n$; let $P=\cup P_i, L=\cup L_i$, $P_i \cap P_j = \emptyset$
 $L_i \cap L_j = \emptyset$
for $i \neq j$ and
let $\lambda$ be a  bijection $\lambda: P\to L$.

A set $\mathcal{K}$ of $k$-tuples $(x_1,x_2, \ldots, x_k)$, $x_i \in P$,
will be called a {\em polygonal presentation} over $P$ compatible
with $\lambda$ if

\begin{itemize}

\item[(1)] $(x_1,x_2,x_3, \ldots ,x_k) \in \mathcal{K}$ implies that
   $(x_2,x_3,\ldots,x_k,x_1) \in \mathcal{K}$;

\item[(2)] given $x_1,x_2 \in P$, then $(x_1,x_2,x_3, \ldots,x_k) \in 
\mathcal{K}$
for some $x_3,\ldots,x_k$ if and only if $x_2$ and $\lambda(x_1)$
are incident in some $G_i$;

\item[(3)] given $x_1,x_2 \in P$, then  $(x_1,x_2,x_3, \ldots ,x_k) \in 
\mathcal{K}$
    for at most one $x_3 \in P$.

\end{itemize}

If there exists such $\mathcal{K}$, we will call $\lambda$ a 
{\em basic bijection}.

Polygonal presentations for $n=1$, $k=3$ were listed in $\cite{Cart}$
with the incidence graph of the finite projective plane of order two 
or three as the graph $G_1$. Some polygonal presentations
for $n>1$ can be found in $\cite{V}$.

\subsection{Construction of polyhedra.}

One can associate  a polyhedron $X$ on $n$ vertices with
each polygonal presentation $\mathcal{K}$ as follows:
for every cyclic $k$-tuple $(x_1,x_2,x_3,\ldots,x_k)$ from
the definition
we take an oriented $k$-gon on the boundary of which
the word $x_1 x_2 x_3\ldots x_k$ is written. To obtain
the polyhedron we identify the sides with the same label of our
polygons, respecting orientation.
We will say that the
polyhedron $X$ {\em corresponds} to the polygonal
presentation $\mathcal{K}$.

The following lemma was proved in \cite{V}:
\begin{lem} \label{Main}
 A polyhedron $X$ which corresponds to
a polygonal presentation $\mathcal{K}$ has
  graphs $G_1, G_2, \ldots, G_n$ as the links.
\end{lem}

\noindent
{\bf Remark.} Consider a  polygonal
presentation $\mathcal{K}$. Let $s_i$ be the number of vertices
of the graph $G_i$ and $t_i$ be the number of edges of $G_i$,
$i=1,...,n$.
If the polyhedron $X$  corresponds to the polygonal
presentation $\mathcal{K}$, then $X$ has $n$ vertices
(the number of vertices of $X$ is equal to the number of graphs),
$k \sum_{i=1}^n s_i$ edges and $\sum_{i=1}^n t_i$ 
faces, all faces are polygons with $k$ sides.

\medskip

\section{Main Construction.}
\subsection{Crucial lemma}
 Let $G$ be a finite projective plane and let $P$ resp. $L$ denote the
set of  its points resp. lines.

   Assume that a bijection $T:P\to L$ is given and satisfies the following
properties 
\begin{enumerate}
\item For each $x\in P$ the point $x$ and the line $T(x)$ are not incident.

\item For each pair $x_1,x_2$ of different points  in $P$ the points
$x_1,x_2 $ and $T(x_1)\cap T(x_2)$ are not collinear.   
\end{enumerate}

\begin{lem} \label{easy}
 Let $T:P\to L$ be as above, $y\in L$ a line. Then the map 
$T^*:I(y) \to I(y)$ given by $T^*(x)=T(x)\cap I(y)$ is a
bijection.  
\end{lem}
\begin{proof}
 By the first property of $T$ the map $T^*$ is well defined, 
by the second property it  must be injective. Since $I(y)$ is
finite, the statment follows.  
\end{proof}

 Let $G,P,L, T:P\to L$ be as above. Let   $P =\{x_1,..x_p \} $ be a labelling 
of points in $P$  and set $y_i=T(x_i)$.
 Consider the following
set  $O\subset P\times P \times P$, consisting of all
triples $(x_i,x_j,x_k)$ satisfying  $x_i\in y_k$, $x_j\in y_i$ and
 and $x_j \in y_k$.

\begin{rem}
 The conditions on $(x_i,x_j,x_k)\in K$ are not cyclic.
We require $x_j \in y_k$ and not $x_k \in y_j$ !!  For this
reason in the polygonal presentations defined below dual 
graphs of $G$ appear. 
\end{rem}

The following lemma is crucial for the later construction:
\begin{lem} \label{crucial}
 A pair $(x_i,x_k)$ resp. $(x_i,x_j)$ resp. $(x_j,x_k)$ is a part
of at most one  triple $(x_i,x_j,x_k)\in K$ and such a triple exists
iff $x_i \in y_k$ resp. $x_j \in y_i$ resp. $x_j \in y_k$ holds.  
\end{lem}

\begin{proof}
 The conditions stated at the end are certainly necessairy. 

1) Let $x_i\in y_k$ be given. Then $y_i$ and $y_k$ are different and the
point $x_j  =y_i\cap y_k$ is uniquely defined. 

2) Let $x_j \in y_i$ be given. Then $x_j$ and $x_i$ are different, so
there is exactly one line $y_k$ containing $x_j$ and $x_l$.

3) Let $x_j \in y_k$ be given. Then $(x_i,x_j,x_k)$ is in $K$ iff
for the map $T^*:I(y_k)\to I(y_k)$  of \lref{easy} 
the equality $T^* (x_i)=x_j $ holds. By \lref{easy}  the point
$x_i$ is uniquely defined. 
\end{proof}

\subsection{Euclidean polyhedra}
 Now we are ready for the polygonal presentations.
Let the notations be as above, $G_1$ and $G_2$ two projective planes with
isomorphisms $J^t:G\to G_t$ and $G_3$ a projective plane with an
isomorphism $J^3:G' \to G_3$ of the dual projective plane $G'$ of $G$.
For $t=1,2$ we set $x^t_i =J^t(x_i)$, $y^t_i=J^t(y_i)$ and for $t=3$
we set $x^3_i =J^3 (y_i)$ and $y^3_i =J^3(x_i)$.

 Let $P_t$ resp. $L_t$ be the set of lines of $G_t$. For
$P=\cup P_t$ and $L=\cup L_t$ we consider the bijection
$\lambda :P\to L$ given by  $\lambda (x_i^t) = y_i^{t+1}$
 ($t+1$ is taken modulo $3$).

 Now consider the subset  $\mathcal{T}$ of $P\times P \times P$ consisting
of all triples $(x^1_i,x^2_j,x^3_k)$ with $(x_i,x_j,x_k)\in K$ and
all cyclic permutation of such triples.

 The stament of \lref{crucial} can be now reformulated as:

\begin{prop}
 The subset  $\mathcal{T}$ of $P\times P \times P$ defines a
polygonal presentation compatible with $\lambda$.
\end{prop}

The polyhedron $X$ which corresponds to 
$\mathcal{T}$ by the construction of \lref{Main} has 
triangular faces and exactly  three vertices with two links
naturally isomorphic to $G$ and one link naturally isomorphic
to the dual $G'$ of $G$.  By \cite{BB} or \cite{CL} the universal 
covering of $X$  is a 
Euclidean building.

\subsection{Hyperbolic polyhedra}

  We continue to use the same notation. We have a projective plane
$G$,  with points $P=\{ x_1,...,x_p\} $ and lines $L =\{y_1,...,y_p \} $
and a subset $K\subset P\times  P \times P$.

Let $w=z_1 \ldots z_n$ be a word of length $n$ in three letters 
$a,b,c$ with $z_1=a,z_2=b,z_3=c$ that does  not contain
proper powers of the letters $a,b,c$. (I.e. $z_z\neq z_{t+1}$ and
$z_n\neq a$). For example $w=abcbcab$ is a possible choice.

 Set $Sign(ab)=Sign(ba)=Sign(ac)=1$ and $Sign(cb)=Sign(ca)=Sign(ba)=-1$.
For $t=1,...,n$ let $G_t$ be isomorphic to $G$ resp. to $G'$ if
$Sign (z_tz_{t+1})=1$ resp. $Sign (z_tz_{t+1})=-1$.

 Fixed isomorphisms induce as above a natural labelling
of the points and lines of $G_t$: $P_t=(x_1^t,....,x_q^t)$ and
$L_t =(y_1^t,....,y_q^t)$.

 For $P=\cup P_t$ and $L=\cup L_t$ we define a basic bijection
$\lambda : P\to L$ by $\lambda (x_i^t)= y_i ^{t+1}$.

 For each triple $(x_i,x_j,x_l)\in K$ we consider
the unique $n$-tuple in $P^n$ such that at the $t$-th place
stands $x_i^t$ resp. $x_j^t$ resp. $x_k^t$ if $z_t$ is equal
to $a$ resp. $b$ resp. $c$. Consider the subset $T_n \in P^n $ of
all such tuples together with all their cyclic permutations.

From \lref{crucial} we immediatly see:

\begin{prop}
 The subset $T_n \in P^n$ is a polygonal presentation over $\lambda$.
By \lref{Main} it defines a polyhedron $X$ whose faces are $n$-gones
and whose $n$-vertices have as links $G$ resp. $G'$.   
\end{prop}

\section{An algebraic construction}

 Let $F=F_q$ be a finite field of charakteristik $p \neq 3$ with $q$
elements. 
Consider the field $K= F_{q^3}$ as  an extension of $F$
of degree $3$. In the sequel we shall denote by $g$ elements of
$K$ and by $a,b,c$ elements of $F$ and call them scalars.
 We denote by $Gr_1$ resp. $Gr_2$ the set of $1$- resp. $2$-dimensional
$F$ vector spaces of $K$.

 The multplicative group $K^*$  operates on the sets $Gr_1$ and
$Gr_2$ by multiplication. The kernel of this operation is precisely
$F^*$ and $K^* / F^*$ operates on both sets simply transitively. 
Especially we can write each element of $Gr_1$ as $gF$ for some
$g\in K^*$.

Let $Tr$ be the trace map $Tr:K\to F$ of the extension $F\subset K$.

 Denote by $E \in Gr_2$ the  $2$-dimensional kernel of $Tr:K\to F$.
We define a map $T:Gr_1 \to Gr_2$ by $T(gF)=gE$. The map $T$ is well defined
bijective and $K^*$ invariant.

\begin{prop}(A.Lytchak, private communication)

For the map $T:Gr_1 \to Gr_2$ and arbitrary $l \neq l_1\in Gr_1$ holds:
\begin{enumerate}

\item  The image $T(l)$ does not contain $l$.

\item  The $l,l_1$ and $T(l) \cap T(l_1)$ generate the vector space $K$.
\end{enumerate}
\end{prop}

\begin{proof}
 Since $T$ is $K^*$ invariant, we may assume $l=F$. Since $Tr(1)=1$, 
$F$ does not lie in $T(F) = E$. Now assume that $l_1 = gF$. 
If the statment is wrong, some non zero element of the form $bg-a$ must
be in $T(F) \cap T(gF) =E\cap gE$. Since $1$ is not in $E$ and $G$ is
not in $gF$, we may assume (replacing $g$ by a scalar multiple ) that
this non zero element is $g-1$. So $g-1 \in E$ and $g-1 \in gE$.

The first inclusion is equivalent to $Tr(g)=1$ and the second one
to $Tr(\frac 1 g)=1$.  Let's prove, that 
if for an element $g\in K^*$ the equalities $Tr(g)=Tr(\frac 1 g) =1$
hold, then $g$ is equal to $1$. 
Assume $g\neq 1$. Then $g$ is not in $F$. Let $m(x)=x^3+ax^2+bx+c$ be the
minimal polynom of $g$. 
Then $c\neq 0$ and $\bar m(x)=x^3 +\frac b c x^2 + \frac a c x +\frac 1 c$
is the minimal polynom of $\frac 1 g$. The condition  
$Tr(g)=Tr(\frac 1 g) =1$ means $a=\frac b c = -1$.
I.e. $m(x)=x^3-x^2 +bx-b= (x^2 +1)(x-b)$ is reducible.
Contradiction. So, $g=1$.

Now we get
a contradiction to $l\neq l_1$. 
\end{proof}

\begin{cor}
 For the  projective plane $\pone$ over finite field
$\fone$ of charakteristique $\neq 3$
there is a bijection $T$ between the set $P$ of 
points and the set $L$ of lines, $T:P\to L$,
that satisfies the following
properties 
\begin{enumerate}
\item For each $x\in P$ the point $x$ and the line $T(x)$ are not incident.

\item For each pair $x_1,x_2$ of different points  in $P$ the points
$x_1,x_2 $ and $T(x_1)\cap T(x_2)$ are not collinear.   
\end{enumerate}
\end{cor}

\end{document}